\input amstex
\documentstyle{amsppt}
\magnification=\magstep1 \hsize=5.2in \vsize=6.8in

\centerline{\bf ON A QUESTION OF D. SHLYAKHTENKO}

\vskip 0.2in \centerline {\rm by} \vskip 0.05in \centerline {\rm
IONUT CHIFAN and ADRIAN IOANA\footnote{The second author was supported by a Clay Research Fellowship}}

\address Math Dept, UCLA, LA, CA 90095-155505
\endaddress

\email ichifan\@math.ucla.edu\endemail
\address Math Dept, UCLA, LA, CA 90095-155505
\endaddress
\email  aioana\@caltech.edu \endemail 
\topmatter \abstract In this
short note we construct two countable, infinite conjugacy class
 groups which admit free, ergodic, probability measure
preserving orbit equivalent actions, but whose group von Neumann
algebras are not (stably) isomorphic. 
\endabstract

\endtopmatter

\document

\head {Introduction.}
\vskip 0.1in
 \endhead Two countable, discrete groups
$\Gamma$ and $\Lambda$ are {\it orbit equivalent} if they admit
free, ergodic, probability measure preserving actions which generate
isomorphic equivalence relations. They are {\it $W^*-$equivalent} if their group von Neumann algebras are isomorphic. D.
Shlyaktenko noticed that, for all known examples of orbit equivalent
groups, their group von Neumann algebras are isomorphic and
speculated that this might be the case in general. Subsequently, a
number of people have also asked the question of whether orbit equivalence
of groups implies $W^*-$equivalence ([Oz06],[Po09]). Motivated by this question, we prove:

\proclaim {Theorem} There exist two countable, discrete, infinite
conjugacy class groups $\Gamma$ and $\Lambda$ which are orbit
equivalent but not $W^*-$equivalent. Moreover, the group von
Neumann algebras $L\Gamma$ and $L\Lambda$ are not stably isomorphic.
\endproclaim
The construction of the groups $\Gamma$ and $\Lambda$ is based on
the observation that being an infinite conjugacy class  group is not
an orbit equivalence invariant. Indeed, by Dye's theorem, $\Gamma_0=S_{\infty}$, the group of finite permutations of $\Bbb N$, is orbit equivalent to $\Lambda_0=\Bbb Z$ ([OW80]).
 This example already
shows that there are orbit equivalent groups which are not $W^*$-equivalent. Further, notice that the map $\Gamma_0\rightarrow \Gamma= (\Gamma_0\times\Bbb F_2)\star\Bbb Z$ turns every pair  $(\Gamma_0,\Lambda_0)$ of orbit equivalent groups into a pair $(\Gamma,\Lambda)$ of orbit
equivalent, infinite conjugacy class groups. 
Finally, by applying the Kurosh type results for free product von Neumann
algebras from [Oz05] we derive that the group von Neumann
algebras of $\Gamma$ and $\Lambda$ are not stably isomorphic.

\vskip 0.2in \head  {Proof of theorem.}\endhead 
\vskip 0.1in
Before proving the theorem, we recall the notion of stable isomorphism of II$_1$ factors. For a II$_1$ factor $M$ and $0<t\in \Bbb R$,
the amplification $M^t$ is defined as the isomorphism class of
$p(\Bbb M_n(\Bbb C)\otimes M)p$, where $n>t$ is an integer
and $p \in \Bbb M_n(\Bbb C)\otimes M$ is a projection of trace $\frac{t}{n}$. It is well known that
this isomorphism class does not depend on the choices of $n$ and $p$.
Then two II$_1$ factors are called {\it stably} isomorphic if
one of them is isomorphic with an amplification of the other.
\vskip 0.1in

{\it Proof.} Let $\Gamma_0$ and $\Lambda_0$ be two infinite amenable groups and assume that $\Gamma_0$ is infinite conjugacy class (ICC) while $\Lambda_0$ is abelian. By [OW80], $\Gamma_0$ and $\Lambda_0$ are orbit equivalent. Further, by [Ga05, Section 2.2], the ICC groups $\Gamma=(\Gamma_0\times\Bbb F_2)\star \Bbb Z$ and $\Lambda=(\Lambda_0\times\Bbb F_2)\star\Bbb Z$ are orbit equivalent. We claim that the group von Neumann algebras $M=L\Gamma$ and $N=L\Lambda$ are not stably isomorphic. Let $M_0=L\Gamma_0$ and $N_0=L\Lambda_0$ and note that $M=(M_0\overline{\otimes}L\Bbb F_2)\star L\Bbb Z$ and $N=(N_0\overline{\otimes}L\Bbb F_2)\star L\Bbb Z$.

If we suppose by contradiction that this is not the case, then we can find  an isomorphism $\theta:M^t\rightarrow N$, for some $t>0$. Since $M_0$ is a factor we can view $M_0^t$ as a subfactor of $M^t$. The commutant of $M_0^t$ in $M^t$ is then equal to $L\Bbb F_2$ ([Po83]). Since the latter is a non-injective factor, by applying [Oz05, Theorem 3.3.] we deduce that there is a unitary $u\in N$ such that $u\theta((M_0\overline{\otimes}L\Bbb F_2)^t)u^*\subset N_0\overline{\otimes}L\Bbb F_2.$ By replacing $\theta$ with Ad$(u)\circ\theta$ we can therefore assume that $\theta((M_0\overline{\otimes}L\Bbb F_2)^t)\subset N_0\overline{\otimes}L\Bbb F_2.$

Since the center of $N_0$ is diffuse, we derive that the commutant of $(M_0\overline{\otimes}L\Bbb F_2)^t$ in $M^t$ is diffuse. However, by [Po83] the commutant of $(M_0\overline{\otimes}L\Bbb F_2)^t$ in $M^t$ is equal to the center of $(M_0\overline{\otimes}L\Bbb F_2)^t$. As $M_0$ is a factor, this gives a contradiction. \hfill$\square$

\vskip 0.1in
{\it Acknowledgment}. We are grateful to Professors Sorin Popa and Yehuda Shalom for useful discussions and encouragement.

\head  References\endhead
\item {[Ga05]} D. Gaboriau: {\it Examples of groups that are measure equivalent to the free group}, Ergodic Theory Dynam.
Systems {\bf 25} (2005), no. 6, 1809–-1827.
\item {[OW80]} D. Ornstein, B. Weiss: {\it Ergodic theory of amenable groups.} I. {\it The Rokhlin lemma.}, Bull. Amer. Math. Soc. (N.S.) {\bf 1} (1980), 161--164.
\item {[Oz05]} N.Ozawa: {\it A Kurosh type theorem for type II$_1$ factors},
Int. Math. Res. Not. (2006), Volume 2006, 1--21, Article ID97560
\item {[Oz06]} N. Ozawa: {\it Amenable Actions And Applications},
Proceeding of the ICM 2006, Vol. II, 1563--1580.
\item {[Po83]} S. Popa: {\it Orthogonal pairs of $\star$-subalgebras in finite von Neumann algebras}, J. Operator Theory {\bf 9} (1983), 253--268.
\item{[Po09]} S. Popa: {\it Revisiting some problems in $W^*-$rigidity}, available at

http://www.math.ucla.edu/~popa/workshop0309/slidesPopa.pdf.

\enddocument